\documentclass[numreferences]{article}

\usepackage{amsgen, amsmath, amsfonts, amsthm, amssymb, enumerate,amscd}
\newtheorem{defn}{Definition}[section]
\newtheorem{prop}[defn]{Proposition}
\newtheorem{lem}[defn]{Lemma}
\newtheorem{thm}[defn]{Theorem}
\newtheorem{cor}[defn]{Corollary}
\newtheorem{rem}[defn]{Remark}

\newcommand {\ZZ}{{\Bbb Z}}
\newcommand {\W}{{\bf W}}
\newcommand {\WW}{{\cal W}}
\newcommand {\B}{{\bf B}}
\newcommand {\X}{{\bf X}}
\newcommand {\XX}{{\cal X}}
\newcommand {\K}{{\cal K}}
\newcommand {\E}{{\cal E}}
\newcommand {\D}{{\cal D}}

\newcommand {\A}{{\cal A}}

\newcommand {\Z}{{\cal Z}}
\newcommand {\C}{{\cal C}}
\newcommand {\Q}{{\Bbb Q}}
\newcommand {\OO}{{\cal O}}
\newcommand {\lprime}{{\ell^2}}
\newcommand {\e}{{\bar e}}
\newcommand {\HH}{{\frak  H}}
\newcommand {\LL}{{\cal L}}
\begin{document}

\title{Relations among Heegner cycles on familes of abelian surfaces}  
\author{Ramesh Sreekantan\\ Institute for Advanced Study \footnote{{\em Current Address:} Department of Mathematics, Duke University} }
\maketitle

\begin{abstract}
We compute relations of rational equivalence among special codimension
$2$ cycles on families of abelian surfaces using  elements of a higher
chow group.
\end{abstract}

\section{Introduction}
 
A conjecture of Beilinson and Bloch asserts that if ${\bf V}$ is a
smooth, geometrically  irreducible, projective variety  defined over a
number field $F$, then
$$   rank     \,  CH^{p}_{hom}({\bf V}_{F})=ord_{s=p}L_F(H^{2p-1}({\bf
  V}),s).$$
Here  $CH^{p}_{hom}$ denotes  the   group of  codimension $p$   cycles
homologically equivalent to   zero modulo rational equivalence  on the
variety.  In particular, this says that the  rank of the Chow group is
finite.  In the   case when $p=1$ the  finiteness  of the rank is  the
Mordell-Weil   theorem.  However very   little  is  known  about  this
conjecture in the case when $p > 1$.  The  purpose of this paper is to
make some progress towards this conjecture in  the case of codimension
2 cycles on certain types of varieties.

The varieties in question are compactifications of modular families of
abelian surfaces over modular  or Shimura curves.  Schoen \cite {sc1},
in the case of modular curves and later Besser \cite {be}, in the case
of  Shimura curves, constructed   certain codimension 2 nullhomologous
cycles supported  in fibres over complex  multiplication points on the
modular curve.  (In fact they showed  that these cycles are in general
not  algebraically equivalent to  0  so they  are  actually in  the so
called  Griffiths  group). These `CM  cycles'   are the  analogues  of
Heegner points on modular curves and are similarly defined over number
fields.   Further, one   can take certain   traces   of these  cycles,
analogous to Heegner divisors, to get cycles defined over fixed number
fields  $F$ and hence  they give  elements of the  groups in question.
The   purpose of this paper    is to construct  relations of  rational
equivalence  between these cycles  in the hope  that one can construct
enough relations to prove that the groups are finite dimensional.

Our  strategy  is to  use the localization   sequence  for higher Chow
groups \cite  {bl1}.  From that one  sees that  to construct relations
between codimension 2 cycles it is enough to construct elements of the
group  $CH^2(\A_{\eta},1)$ where $\A_{\eta}$  is  the  generic abelian
surface of  the family and then compute  their boundary.  It turns out
that there are some very natural elements constructed by Collino \cite
{co}  and interestingly, their boundary  can be expressed  in terms of
the CM cycles.

We then have the following:
\begin{thm}

Let ${\bf X(D_0,N)}$ be a Shimura curve parametrising abelian surfaces
with endomorphism ring an Eichler order (loc. cit. Section 2) of level
$N$ in a division algebra of discriminant $D_0$.  Let ${\bf W(D_0,N)}$
denote  (the  non-singular compactification of)  the universal abelian
surface  over ${\bf X(D_0,N)}$.  Let  $p,a,b$ be  the invariants which
determine Hashi -moto's model (ibid.)  and let $P$  and $Q$ denote two
$2$-torsion points in the generic fibre.

Then there  are relations in  $CH^2({\bf W(D_0,N)})\otimes \Q$  of the
form
$$  \sum_{[\LL]} \epsilon_{P,Q}([\LL])   \sum_{r,s}  \sum_{d^2|\Delta}
d \Z_{\Delta/d^2,[\LL],s_1/d,-s_1/d} \equiv 0$$
where $\LL$ runs through all even theta characteristics, $r,s \in \ZZ$
such that $\frac{n^2p-r^2}{4} \in \ZZ_{\geq 0}$, $\frac{4bDn^2-s^2}{4}
\in \ZZ_{\geq 0}$, $n$ is coprime with $2k$, $r,s,n$ mutually coprime,
$s_1=p(s/2)-aDr)$ and
$$\Delta=\frac{ps^2-4aDrs+4bDr^2-4Dn^2}{4}$$
with  $d$ running through all $d^2|\Delta$   such that $\Delta/d^2$ is
still a discriminant.

\end{thm}
Here $\Z_{\Delta}$ is a {\em Heegner Cycle of discriminant $\Delta$}
(loc.  cit Section 5  ), a codimension  2 cycle which is homologous to
zero in  $CH^2({\bf W(D_0,N)})\otimes \Q$  and $\epsilon_{P,Q}([\LL])$
is a sign function which  depends on the level and   the points $P$  and
$Q$ (loc. cit. Section 4).

Relations  between Heegner {\em   points} on modular curves  have been
studied before by many other authors. In Zagier \cite {za} and van der
Geer \cite{vdg}, they construct relations  between these points  which
are very similar to ours, the main difference being that our relations
involve the level 2 structure while theirs do not.

As is known   \cite{za},\cite{gz},  Heegner points are   very  closely
related  to coefficients of modular  forms of weight $\frac{3}{2}$ and
relations  of  rational   equivalence  between Heegner  points   imply
relations between  coefficients  of these  forms.   Similarly, Heegner
cycles are supposed to be related to  coefficients of modular forms of
weight $\frac{5}{2}$.  Some evidence of this can be found in
\cite{sc2},\cite{zh} and \cite{ne}. So  relations of rational 
equivalence between these cycles should give rise to relations between
coefficients of such modular forms.

More   interestingly, recently  Borcherds  \cite{bo} has  used his new
constructions  of  automophic  forms to  construct  relations  between
Heegner points on  modular   and Shimura curves  and   more generally,
relations between special divisors  on   modular varieties, and   {\em
prove} that these special  subvarieties are related to coefficients of
modular forms.   Starting  from certain  meromorphic modular  forms he
constructs other automorphic forms whose divisors  are on such special
subvarieties.  It would be very interesting  to see if his methods can
be generalized to construct  Collino type elements  of the higher Chow
groups and use  them to prove that the  Heegner cycles are related  to
coefficients of modular forms. At the moment it is not clear to me how
this can be  done, but a  weaker statement, namely computation  of the
regulator of Collino's element may be a little more tractable.

The outline of the paper  is as follows:  In Section 2 we describe the
varieties and cycles in question. In Section 3 we introduce the higher
Chow groups and the localization sequence that we use. In Section 4 we
describe Collino's construction  of  the elements of  the  higher Chow
groups.  In  Section  5 we compute the  boundary  of the elements  and
relate them to the  CM cycles. Finally, in  section 6 we put them  all
together to get our main result and give some examples.

{\em Acknowledgements}: I would like to thank my advisor Spencer Bloch
for suggesting  the problem, Madhav Nori for  refering me to Collino's
work, Patrick  Brosnan  for   helping  me  get started  and  Najmuddin
Fakruddin   for his invaluable   help  and innumerable suggestions and
corrections.  I  would  also  like  to thank  the  hospitality of  the
Institute for Advanced Study where this  paper was completed. Finally,
I  would like to   thank Chad Schoen  and the  referee  for their very
useful comments.

\section{Preliminaries}

\subsection{Eichler Orders}

The varieties we consider  are universal families of  abelian surfaces
whose endomorphism  rings contain an  order in a quaternionic division
algebra.  To describe precisely we use the theorems of Hashimoto
\cite{ha} and Roberts \cite{ro}. Another standard reference is \cite{vi}.

Let ${\bf B}$ be an indefinite quaternion algebra over ${\Bbb Q}$ with
discriminant $D_0$. An {\bf Eichler Order of level $D=D_0N$} is an order
$\OO$ in $\B$ such that
$$\OO\otimes    \ZZ_{\ell}  \equiv  \pmatrix   *& *\\0 & * \endpmatrix
\,\,mod\,\,\, \ell\,\,\,\,\, for \,\,\, all \,\,\, \ell|N$$

Let ${\cal O}$ be an  Eichler order of level $D=D_0N$  where $N$ is  a
positive integer prime to $D_0$.  Let ${\cal S}$  be the set of primes
where the division algebra is ramified  ( i.e those primes $\ell$ such
that  ${\bf B} \otimes {\Bbb Q}_{\ell}   \neq M_2({\Bbb Q}_{\ell}) $).
Choose   an   auxiliary  prime $p$    such  that the    Hilbert symbol
$(-D,p)_{\ell}=-1$ if  and only if $\ell\in {\cal  S}$.  Such a $p$ is
guaranteed  by    Dirichlet's  theorem   on  primes     in  arithmetic
progressions.  Let $a$ and $b$ be such that $a^2D+1=bp$.  Then we have
the following theorem described  in Hashimoto, \cite{ha}, though parts
of it existed in some form in earlier papers.

\begin{thm}

Let   $\B$    be a quaternionic     division  algebra  of discriminant
$D_0$. Then: 

1. $\B \simeq \Q \oplus  \Q i \oplus \Q j \oplus \Q ij$ where
$i^2=-D$, $j^2=p$ and $ij=-ji$.

2. The order ${\cal O} \simeq \ZZ e_{1} \oplus \ZZ e_2 \oplus \ZZ e_3
\oplus    \ZZ e_4$  where    
$e_1=1,e_2=(1+j)/2, e_3=(aDj+ij)/p$ and $e_4=(i+ij)/2$

3. There is a skew-symmetric pairing on  $\B$ which is $\ZZ$ valued on
${\cal O}$  given by $E(x,y)=tr( x  {\bar y} i^{-1})$ where ${\bar y}$
denotes conjugation.  Further, the elements $\eta_1=e_4-\frac{p-1}{2},
\eta_2=-aDe_1-e_3, \eta_3=e_1$ and $\eta_4=e_2$ are a symplectic basis
for ${\cal O}$

\end{thm}

Since $\B$ is indefinite,   any two Eichler   orders of level  $D$ are
conjugate \cite{bd}, so there is no loss of generality in working with
this  model and we  can  fix the isomorphism  of  $\B \otimes {\Bbb R}
\simeq M_2({\Bbb R})$ by defining
$$\Phi_{\infty}(i)=
\pmatrix 0 & -1 \\
         D & 0 \endpmatrix,\,\,\,\,\, \Phi_{\infty}(j)=\pmatrix {\sqrt
         p}&0\\ 0& -{\sqrt p} \endpmatrix $$

We  will  call   this description   of   the quaternion  algebra  {\bf
Hashimoto's model}.

\bigskip
\noindent {\em Orientations}
\bigskip

The above   description  is  not     rigid  enough.  So,     following
\cite{ro},\cite{bd} we rigidify the  definition further by defining an
{\bf Oriented Eichler  Order of level  $D_0N$}  as follows. For  each
prime $\ell$ dividing $D_0$ there are two algebra homomorphisms,
$$ {\frak   o}_{\ell}:\OO  \otimes {\Bbb  F}_{\ell} \rightarrow  {\Bbb
F}_{\ell}$$
$$ \pmatrix a & b \\c & d \endpmatrix \rightarrow a \,\,\,\pmatrix a & b \\c & d \endpmatrix \rightarrow d $$
and similarly  for a prime $\ell$ dividing  $N$ there are two distinct
algebra homomorphisms
$${\frak   o}_{\ell}:\OO   \otimes {\Bbb   F}_{\ell} \rightarrow {\Bbb
F}_{\lprime} $$
$$ \pmatrix a & b \\ \ell b^{\sigma} & a^{\sigma}\endpmatrix \rightarrow a \,\,\, \pmatrix a & b \\ \ell b^{\sigma} & a^{\sigma}\endpmatrix \rightarrow a^{\sigma}$$
where $\sigma$ is  the non-trivial automorphism of ${\Bbb F}_{\lprime}
$. An oriented Eichler Order is an Eichler order along with a choice of
one  of  these homomorphism ${\frak  o}_{\ell}$,  called a $\ell$-{\em
orientation }, for each $\ell$ dividing $D=D_0N$.

\subsection {Shimura Curves}

With the model of  $\B$,$\OO$, the isomorphism $\Phi_{\infty}$ and the
orientation we can  form the   family of  surfaces we are   interested
in. Let $\Gamma(1)=\Gamma_{\OO}(1)$ be the group
$$\Gamma_{\OO}(1)=\{x \in \OO|Nm(x)=1\}$$
where  $Nm$ denotes the reduced norm.    $\Gamma(1)$ acts on the upper
half plane ${\frak  H}$ through the   embedding $\Phi_{\infty}$.  The
quotient  is an algebraic  curve   whose complex points represents   a
component of  the moduli of abelian  surfaces whose endomorphism rings
contain  $\OO$. In   this  case   the  generic  abelian surface    has
endomorphism ring actually equal to $\OO$.

We  further assume  that the abelian   surfaces have  full level  $2k$
structure for  some  $k>1$. This is because  we  will use  the level 2
structure and we need further level structure to  ensure that there is
a universal  family over the curve.  The $k$-level structure  will not
play any other part in our calculations, so will be suppressed but not
forgotten in the remaining part of this paper.

The universal family can be described as follows. Let 
$$\Gamma(2k)=\Gamma_{\OO}(2k)=\{\gamma  \in    \Gamma(1)|\gamma    \in
Ker(\Gamma(1) \rightarrow Aut(\OO/2k\OO))\}$$ 
The  {\em   Jacobi  Group}  $\Gamma(2k)    \ltimes  {\cal  O}$    with
multiplication defined     by   $(\gamma,x)   \star    (\mu,y)=(\gamma
\mu,x\mu+y)$ acts on ${\frak H} \times {\Bbb C}^2$ via
$$(\pmatrix         a&b        \\         c&b\endpmatrix,x)      \cdot
(\tau,z_1,z_2)=(\frac{a\tau+b}{c\tau+d},       (c\tau+d)^{-1}(\pmatrix
z_1\\z_2
\endpmatrix + \Phi_{\infty}(x) \pmatrix \tau\\1 \endpmatrix )) $$
and  the quotient is a  complex threefold  which  is a family of abelian
surfaces over the Shimura curve ${\frak H}/ \Gamma(2k)$.

The fibre over a point $\tau$ is 
$$\A_{\tau} = {\Bbb  C}^2/ \Phi_{\infty}({\cal O})\pmatrix \tau \\  1
\endpmatrix$$ 
and $\A_{\tau} \simeq \A_{\tau'}$ if  and only if $\tau=\gamma(\tau')$
for  some $\gamma  \in  \Gamma(1)$,   with their level   structures
coinciding if and only if $\gamma \in \Gamma(2k)$.

Let ${\bf Y}={\bf  Y(D_0,N)}$ denote the curve ${\frak H}/\Gamma(2k)$.
If $D_0=1$ the curve  is  not compact but   it can be  compactified by
adding finitely  many cusps.   Let   $\X={\bf X(D_0,N)}$  denote   the
compactified  curve.   Let     $\W={\bf    W(D_0,N)}$  denote      the
desigularisation  of  the self product  of   the universal family over
${\bf   X(D_0,N)}$    as     described     in     \cite{sc1}.      Let
$\WW=\WW(D_0,N)=\A_{\eta}$ denote   the generic fibre  of this family.
As remarked earlier, it is  an abelian surface with endomorphism  ring
$\OO$.

\bigskip
\noindent {\em Polarizations}
\bigskip

There    is      an  involution    on  ${\cal      O}$     defined  by
$d^{\dagger}=i^{-1}{\bar d} i$ where $i$  is as in Hashimoto's theorem.
The skew symmetric pairing on ${\Bbb C}^2$ defined by
$$<\Phi_{\infty}(x)\pmatrix \tau \\ 1 \endpmatrix,
\Phi_{\infty}(y)\pmatrix \tau \\1\endpmatrix>=E(x,y) \eqno(1)$$
where  $E$ is the form  coming  from part $3$   of  theorem $2.1$, is  a
Riemann form   for $A_{\tau}$. This  is a  principal  polarization and
$d^{\dagger}$ is  the Rosati involution  corresponding  to it.  A well
known theorem \cite{mu} pg  190, says that the  fixed points under this
involution correspond to  elements of the  N\'eron-Severi group of the
abelian surface.  It is easy to see that  the elements $1$,$j$ and $k$
are  invariant under  this    involution  and  are   clearly  linearly
independent, so from  that one can see  that  the rank  of the generic
N\'eron-Severi group is  $3$. A little more  work gives the following,
which is a simple corollary to Hashimoto's theorem, but will be useful
later on.

\begin{cor}

The elements   $e_1=1$, $e_2=\frac{1+j}{2}$ and $e_3=\frac{aDj+ij}{p}$
give a  basis for the N\'eron-Severi  of the  generic fibre.  Further,
for ample $g$, the cup  product is given by $<g,g>=2Nm(\phi_g)$, where
$\phi_g$  is  the endomorphism corresponding   to $g$. With  that, the
intersection matrix is
$$
\begin{array}{c|ccc}

  & e_1 &e_2 & e_3 \\ \hline e_1 &2 & 1 & 0 \\ e_2& 1& \frac{1-p}{2} &
  aD\\ e_3 & 0 &aD & 2bD

\end{array}
$$
where $a^2D+1=bp$ 
\end{cor}

\begin {proof}

From the Riemann-Roch theorem, one sees that, for an ample divisor $g$,
$$<g,g>=2\chi({\cal      L}(g))    \,\,\,    and   \,\,\,   \chi({\cal
L}(g))^2=deg(\phi_g)$$
where  ${\cal L}(g)$ is the  corresponding line bundle and $\phi_g$ is
the corresponding endomorphism.

This, coupled with the formula for the degree in terms of the norm
$$deg (\phi)=Nm(\phi)^2$$
where $Nm$ is the reduced norm, shows that 
$$<g,g>=\pm2Nm(\phi_g)$$
The sign can be   determined from the  fact  that $<1,1>=2$ where  $1$
represents the class of the principal polarization.

\end{proof}

From now on we pick curves in the generic  fibre which represent these
classes $e_i$ and  we  denote them by   $e_i$ as well. This  will  not
affect  anything   as all  of   our results  will   be modulo divisors
homologous to zero in the  fibres. We also fix  the polarization to be
the one in equation (1) and  denote it by $\Theta$.

\subsection { Humbert Invariants}

A lemma in \cite{be2}, Lemma 4.1 shows that there is an equivalence of
categories between the   category of even lattices  of  rank n with  a
element $\zeta$ of norm 2 with the category  of lattices of rank (n-1)
such that the quadratic form represents only numbers $\equiv 0$ or $1$
mod  $4$.  The correspondence   is   given  by  
$$(M,\zeta) \rightarrow N(M,\zeta) \hspace{1in} N \rightarrow (M_N,\zeta)$$
where 

$\bullet$   $N(M,\zeta)$  is the lattice  such  that   $2N(-2)$ is the
orthogonal complement to $\zeta$ in $2M+\ZZ\zeta$
 
$\bullet$ $M_N$ is  the sublattice of $(N(-2) \oplus  \ZZ\zeta)\otimes
\Q$  containing  $N\oplus \ZZ\zeta$  and  all   elements  of the  form
$\frac{1}{2}(x+(x,x)\zeta)$ with $x$ in $N$.

Here $N(k)$ denotes  the lattice which  has the same underlying  $\ZZ$
module as $N$ but with the pairing multiplied by $k$.

The N\'eron-Severi lattice of an  abelian surface is an even  lattice.
A principal   polarization is an element of   norm 2.  So a  choice of
principal polarization  allows  us  to work  with a   smaller lattice,
namely the primitive N\'eron-Severi lattice.
$$NS(\A,\Theta)=\{ v \in NS(\A)|<v,\Theta>=0 \}$$
with the pairing given by
$$(v,v)_{\Theta}=-2<v,v>$$
More generally, one  can work with   the entire $H^2(\A.\ZZ)$ and  the
primitive cohomology, $H^2(\A,\Theta)$.

Classically, this  was first studied by Humbert  \cite{hu}, and  so we
call the  pairing $(,)_{\Theta}$ the {\bf Humbert  norm} and if $v \in
NS(\A,\Theta)$   then   $(v,v)_{\Theta}=(v,v)$ is    the  {\bf Humbert
Invariant of v, H(v)}.

From now on  unless otherwise mentioned, we will  always work with the
Humbert norm  on  a   lattice.  From  the  corollary  above    and the
description of  the   correspondence above we   see that  the  Humbert
lattice for the generic N\'eron-Severi is generated by
$$\e_2:=\frac{1}{2}(2e_2-<e_2,e_1>e_1) \,\,\,\,\,and \,\,\,\,\,  \e_3:
= \frac{1}{2}(2e_3-<e_3,e_1>e_1)$$ 
with intersection matrix $$
\begin{array}{c|cc}

   & {\bar e_2} & {\bar e_3}\\ \hline
    {\bar e_2} & p & 2aD \\
    {\bar e_3} & 2aD & 4bD 

\end{array}
$$
This matrix has determinant 4D.

Humbert studied the   moduli of abelian surfaces whose  N\'eron-Severi
contains an element  with a non-zero Humbert  invariant, which are now
called  {\bf Humbert surfaces}.  This amounts  to saying that the rank
of the  N\'eron-Severi is at  least 2 and   he realised that  this was
equivalent to saying that the endomorphism ring contains an order in a
real   quadratic field.  One   can  easily see  that,  in general, the
endomorphism algebra is simply  the  Clifford algebra of  the  Humbert
lattice.

Humbert's definitions were more analytic and  quite different from the
above   defintions,  which are due  to   Kani.  For some computations,
however, it is more useful to use the analytic definition, so for that
reason we will  give it here. The equivalence  of the two can be found
in Kani \cite{ka}.

Let  $V'$  be  the  set   of integral,  skew   symmetric 4x4 matrices,
considered  as  a subgroup of   $M_4(\ZZ)$.  Let $J$   be the subgroup
generated  by $\pmatrix   0   & I_2  \\-I_2  &0  \endpmatrix$  and let
$V$=$V'/J$.

There is an action of $Sp_4(\ZZ)$ on $V'$ given by
$$ v \rightarrow M^t v M$$
for  $M$ in $Sp_4(\ZZ)$ which leaves  $J$ invariant and hence descends
to an action on $V$. A $v$ in $V'$ looks like
$$v = \pmatrix 0 & a & b_1 & b_2 \\ -a & 0 & b_3 & b_4\\ -b_1 & -b_3 &
0 & d \\-b_2 &-b_4 & -d & 0 \endpmatrix = \pmatrix A & B\\ -B^t & D
\endpmatrix$$
and one  defines  the {\bf Humbert   Invariant} to  be  $$ (v,v)=H(v)=
(b_1+b_4)^2-4(b_1b_4-b_3b_2-ad)=Tr(B)^2-4det(B)+ad$$
$H(v)$ depends  only on  $v$  mod $J$   and  $v$ is  said  to be  {\em
primitive} if  it is not  divisable  by a natural   number $\geq 2$. A
lemma due to Humbert says that two primitive  elements $v_1$ and $v_2$
are equivalent if and only if $H(v_1)=H(v_2)$.

Let ${\bf \tau}=\pmatrix \tau_1 & \tau_2 \\ \tau_2 & \tau_3
\endpmatrix $ be a point in the Siegel upper half space ${\frak H}_2$.  
Every primitive $v$ in $V$ determines a subvariety given by ${\bf
\tau} \in {\frak H}_2$ such that
$$\pmatrix {\bf  \tau} & I_2 \endpmatrix  v \pmatrix {\bf \tau} \\ I_2
\endpmatrix =0 $$
$$\tau A \tau + \tau B -B^t \tau +D=0$$
and if ${\bf \tau}$ satisfies such an equation it is said to satisfy a
{\bf singular relation of invariant H(v)} and the subvariety is called
the Humbert  surface of invariant $H(v)$. $V$  is  called the space of
singular relations. Equivalently, a singular relation can be described
as a relation between the coefficients of $\tau$ of the form
$$v=\alpha\tau_1         +         \beta\tau_2+      \gamma\tau_3      +
\delta(\tau_2^2-\tau_1\tau_3) + \epsilon=0$$
with invariant $H(v)=\beta^2-4\alpha\gamma-4\delta\epsilon$. Sometimes
it  is more convenient to  work with this  description. We will denote
$v$ by  $(\alpha,\beta,\gamma,\delta,\epsilon)$.  In fact  \cite {ka},
Prop.5.4, there is a basis for the primitive part of $H^2(\A,\ZZ)$ for
which    this vector gives  the  corresponding    2-form in the  other
description of the Humbert surface.

\subsection{CM points and CM cycles}

The  cycles we  are interested  are    codimension 2 cycles on   these
families of abelian surfaces.  There  are two basic types,  horizontal
cycles  and vertical cycles.  Horizontal  cycles are roughly images of
the Shimura curve under some section.  Vertical cycles are those which
are supported in fibres  over some points.  In  most fibres,  all such
cycles will come  by restriction of a  divisor on universal  family to
the special fibre.   However, there are some   points in whose  fibres
there are extra elements  in the N\'eron-Severi   of the fibre,  hence
extra codimension 2 cycles in the family.

\bigskip

\noindent {\em CM Points}

\bigskip

As we have seen before, the  rank of the  generic N\'eron-Severi is 3,
so when there is an extra cycle, the rank goes up to 4. If $\tau$ is a
point in the moduli where the rank  of the N\'eron-Severi of the fibre
$\A_{\tau}$ is 4, then  the abelian surface $\A_{\tau}$ is necessarily
isogenous to a product  of elliptic curves with complex multiplication
by an the imaginary quadratic field $K=\Q(\tau)$ and $End_0(\A_{\tau})
\simeq M_2(K)$. Such a point $\tau$ is called a {\bf CM point}.

Determining a CM point is equivalent to determining an embedding of an
imaginary  quadratic   field $q:K \hookrightarrow   \B$  such that $\B
\otimes K \simeq M_2(K)$  as from Prop. 9.4 of  \cite {sh} there is  a
unique point $\tau$ in ${\frak H}$ such that
$$\Phi_{\infty}(q(K^{\times}))=\{\gamma    \in  \Phi_{\infty}(B^{\times}
\cap GL_2^{+}(\Q))\,|\, \gamma\tau=\tau\}$$

We  can normalize  these embeddings  as in  \cite  {sh} (4.4.3).  Once
having normalized the embedding of $K \hookrightarrow \B$ one has also
normalized the embedding  of $K \hookrightarrow M_2(K)=\B \otimes  K$.
$K \cap End(\A_{\tau})$ is an order of discriminant $\Delta$ and hence
will be denoted  by $\OO_{\Delta}$. Through the  embedding $q$ and the
maps ${\frak  o}_{\ell}$ for $\ell$ dividing  $D$, we have also chosen
an  orientation  of the    order  $\OO_{\Delta}$.

\bigskip

\noindent {\em CM cycles}

\bigskip

Let  $\alpha$ be a  traceless   element  of ${\cal O}_{\Delta}$   with
$\alpha^2=\Delta$, so $\alpha$     is purely  imaginary. The    Rosati
involution extends to an involution of $\B \otimes K$ and acts by
$$(d \otimes \alpha)^{\dagger}=d^{\dagger} \otimes \bar
{\alpha}=-d^{\dagger}\otimes  \alpha$$
Observe  that the element  $i \otimes \alpha$   is fixed by the Rosati
involution. We define a {\bf CM cycle class} to be
$$\Z_{\tau}=i\otimes \alpha$$
the class  of this element in the   N\'eron-Severi of $\A_{\tau}$.  By
construction $\Z_{\tau}$  is {\em not}  a generic  class.  Further, in
this case one can see that the intersection form on the N\'eron-Severi
lattice  being  thought  of  as Rosati   fixed  elements  is  given  by
$<v,v>=2det(\phi_{v})$.  From  that one can see that   the CM cycle is
orthogonal   to  the generic  N\'eron-Severi   and the Humbert norm is
$$(\Z_{\tau},\Z_{\tau})=-4D\Delta=-4Ddisc(\OO_{\Delta})$$

\begin {rem} This definition of the CM cycle almost agrees with the 
definitions in \cite{sc1} and \cite{be}, where it is defined to be the
minimal  generator  of the  orthogonal    complement  of the   generic
N\'eron-Severi. The definitions are the same  when $\Delta$ is odd, but
when $\Delta$ is even, our cycle is twice the generator.
\end{rem}

The CM  cycles are the analogues of  CM points  on modular and Shimura
curves and  have many similar properties.  They  are also defined over
number fields.  Schoen, in the modular case and Besser, in the Shimura
curve case, show that they are homologous to zero in $CH^2(\W)$ and so
give rise to  interesting  elements in  $CH^2_{hom}(\W)$.  In fact  he
shows that in general they are non-trivial in the Griffiths group.

\section {Higher Chow groups and the Localization sequence}

In  this section we introduce  the higher Chow  groups that we use and
the Localization theorem.  No proofs are  given.  The proofs of  these
results may be found in \cite{bl1}.

\subsection{The group $CH^2(\A,1)$}

The group  we use  is  $CH^2(\XX,1)$ where $\XX$   is surface.   It is
defined  as  a a  certain  subquotient of the  group  of codimension 2
cycles on $\XX \times {\Bbb A}^{1}$.  A theorem of Bloch's \cite {bl1}
identifies it with   the  ${\cal K}$-cohomology   group $H^1(\XX,{\cal
K}_2)$, where $\K_2$ is the {\em sheaf} coming from the presheaf given
by
$$U \longrightarrow K_2(U)$$
With this description  and the Bloch-Gersten-Quillen resolution of the
sheaf   $\K_2$ one  can  see  that   an  element of  $CH^2(\XX,1)$  is
represented by a formal sum
$$\sum (\C_i,f_i)$$
where $\C_i$ are curves on $\XX$ and  $f_i$ are functions on the $C_i$
such that they satisfy the cocycle condition
$$\sum div(f_i)=0$$

In  $CH^2(\XX,1)$ there are certain elements   coming from the product
structure on Chow groups,
$$CH^1(\XX,1) \otimes CH^1(\XX,0) \rightarrow CH^2(\XX,1) \eqno(1)$$
A theorem of  Bloch's \cite {bl1} shows  that 
$$CH^1      (\XX,1)     =    \Gamma  (\XX,\OO_{\XX}^{*})  \hspace{1in}
CH^1(\XX,0)=CH^1(\XX)=Pic(\XX)$$
These elements are of the form
$$\sum (\C_i,f_i)$$
where $\C_i$ are divisors on  $\XX$  and $f_i$ are constant  functions
hence     automatically    satisfy   the    cocycle    condition    as
$div(f_i)=0$.  Such elements   are  called  {\bf decomposable}.   More
generally one can construct more such elements  by looking at norms of
elements of extensions of    the base field. Let   $CH^2(\XX,1)_{dec}$
denote the subgroup  of $CH^2(\XX,1)$  generated  by such cycles.  The
group
$$CH^2(\XX,1)_{indec}=CH^2(\XX,1)/CH^2(\XX,1)_{dec}$$
is called the subgroup of {\bf indecomposable} cycles and is sometimes
non-trivial.

\subsection {The Localization sequence}

The localization  sequence for higher Chow  groups is  an extension of
the usual localization sequence for  Chow groups to the left.  Suppose
$\XX$  is a smooth projective  variety  and ${\cal Y}  \hookrightarrow
\XX$ is a closed divisor, then we have an exact sequence
$$  ...\rightarrow  CH^2(\XX,1) \rightarrow CH^2(\XX  \backslash {\cal
Y},1) \stackrel{\partial}{\rightarrow}  CH^1({\cal   Y},0) \rightarrow
CH^2(\XX,0) \rightarrow...$$
In particular, if $\XX$ is a family of surfaces over  a curve $\Z$ and
${\cal Y}_{z}$  is the fibre  over a point $z  \in \Z$, then we have a
sequence
$$...\rightarrow CH^2(\XX,1) \rightarrow CH^2({\cal Y}_{\eta},1)
\stackrel{\partial}{\rightarrow}      \bigoplus_{z  \in \Z} CH^1({\cal
Y}_{z}) \rightarrow CH^2(\XX) \rightarrow ..$$
where ${\cal Y}_{\eta}$ is the  generic fibre of  the family. So, from
this our strategy  becomes clear: to  construct relations of rational
equivalence between codimension  2  cycles supported  in fibres  of  a
family of  surfaces it  suffices  to construct elements of  the higher
Chow group of the generic fibre and compute the boundary under the map
$\partial$.   If $\sum   (\C_i,f_i)$  is   an  element  of  $CH^2({\cal
Y}_{\eta},1)$, so now the $\C_i$ are curves in  the generic fibre, the
boundary of the element is
$$\partial(\sum (\C_i,f_i))=\sum div(f_i)$$
where  $div(f_i)$ denotes the  divisor of $f_i$  being though  of as a
function on the closure ${\bar \C_i}$ of $\C_i$. This divisor may have
vertical componenents.  For  example,  if  $(\C,f)$ is  a  decomposable
element in $CH^2({\cal  Y}_{\eta},1)$, where $f$ is  now a  function on
the base, $$\partial( (\C,f) )=\sum a_i\C_{x_i}$$
where $div(f)=\sum a_ix_i$  and $\C_{x_i}$ denotes the restriction  of
$\C$ to the fibre over the  point $x_i$. Note that  in this manner one
can only get  relations between cycles obtained  by the restriction of
generic cycles   to the various    fibres. Since  the cycles   we  are
interested  in are   not  of  this  type,  one   cannot hope  to   use
decomposable   elements for  our  purposes.  However,   on the generic
abelian   surface there exist   indecomposable  elements  In the  next
section we describe these elements constructed by Collino \cite {co}.

\section { Collino's elements of $CH^2(\A,1)$}

\subsection{Construction of the Element}

In  \cite {co}, Collino constructs   certain elements of  $CH^g(\A,1)$
where $\A$ is the  Jacobian of a hyperelliptic   genus $g$ curve.   In
particular,   if $\A$  is a   simple,   principally polarized  abelian
surface, it is the Jacobian of a smooth genus 2 curve,  so one can use
his construction. It is as follows:

Let $\A=Jac(\C)$ where $\C$ is  a genus 2  curve which represents  the
principal polarization.  Since  $\C$    is hyperelliptic there   is  a
function  $f:\C\rightarrow {\Bbb P}^{1}$   such that $div f=2(P)-2(Q)$
where $P$ and $Q$ are two ramification points.

There are  maps    $i_{P}(x)=x-(P)$  and $i_{Q}(x)=x-(Q)$   from   $\C
\hookrightarrow Jac(\C)$. Let  $\C_{P}$ and $\C_{Q}$ denote the images
of  $\C$  under $i_{P}$ and $i_{Q}$  respectively   and let $f_P$ and
$f_Q$ denote  the  function  $f$  being thought  of as    functions on
$\C_{P}$ and $\C_{Q}$ respectively. Then
$$(C_P,f_P)+(C_Q,f_Q)$$
is an element of $CH^2(\A,1)$.  This is because 
$$div(f_P)=2(O)-2(Q-P) \,\,\, and \,\,\,\, div(f_Q)=2(P-Q)-2(O)$$ 
but $(P-Q)=(Q-P)$ in $Jac(\C)$, so
$$div(f_P)+div(f_Q)=2(O)-2(Q-P)+2(P-Q)-2(O)=0.$$ 

Collino   proves that this  element  is  indecomposable on the generic
abelian surface. Note that the exact choice of the function $f$ is not
that important as  any other choice  of function with the same divisor
would give the same class in the group of indecomposable cycles.

\subsection{Computation of the Boundary}

Let $\X$ as before be the Shimura curve with full level $2k$ structure
and $\W$  and $\WW$  be the  universal  family and  the generic  fibre
respectively. Since we have full level 2 structure, by making a choice
of two 2-torsion sections  $P$ and $Q$, we  can apply the construction
in  section    $4.1$  to the   generic fibre   to  get   an element of
$CH^2(\WW,1) =  CH^2(\A_{\eta},1) $.  In  this section we will compute
the boundary of this element.

The boundary of this element of $CH^2(\WW,1)$ is of the form 
$$\sum_{x} a_{x}{\cal D}_{x}$$
where $\D_{x}$ are  codimension 1 cycles in the  fibre over the  point
$x$.  In many fibres, $\D_x$ is simply  the restriction of the generic
curves $\C_{P}$ and $\C_Q$ to the fibre, but in some  cases it is more
interesting. There are two such cases \cite{we}. The first is when the
point is a cusp,  so the fibre is   a degenerate abelian surface.   In
this  case  a   theorem  of    Schoen  \cite{sc1}, essentially     the
Manin-Drinfeld principle, shows   that  these cycles are  torsion   in
$CH^2(\W)$ so one can neglect  them as we  are only interested in  the
rational Chow group.  The second case is when the curve breaks up into
a sum  of two elliptic  curves intersecting at  a point, such that the
two ramification points lie   on different components.  Later on,   in
sections $4.3$     and $5$, we  will   describe   precisely  when this
happens. In  this  case the  cycle is  not simply the  restriction and
involves the   components  in a non-trivial   manner. To  compute  the
boundary we do the following local computation.

\begin{thm}

Let $\tau$ be a point on $\X$ in whose fibre  the genus 2 curve $\C_P$
degenerates  into a sum of two  elliptic  curves $\E_1$ and $\E_2$ and
such that $O$ lies on $\E_1$ and $Q-P$ lies on $\E_2$. Then
$$\partial ((C_P,f_P)+(\C_Q,f_Q))=2(\E_2-\E_1)$$
up to  the boundaries of   decomposable elements and  images of cycles
homologous to 0 in the fibre.

\end{thm}

\begin{proof}

In a neighborhood of the point  $\tau$ the boundary of $(\C_P,f_P)$ is
of the form
$$\partial((\C_P,f_P))=a\E_1+b\E_2+{\cal H}$$
where $a$ and $b$ are in  $\ZZ$ and ${\cal  H}$ denotes the closure of
the horizontal sections of $div(f_P)$.  Locally,  we can always find a
decomposable element  whose boundary is  of the form $b(\E_1+\E_2)$ so
subtracting this element from Collino's  element  allows us to  assume
that the boundary is of the form
$$\partial((\C_P,f_P))=a\E_1+{\cal H}$$
To show that the boundary is not simply the restriction of the closure
of a divisor in the generic fibre amounts  to showing that $a \neq 0$.
To do this we intersect  with $\E_2$ and  use the fact that a function
restricts  to a divisor  of degree 0 in  a  cycle not contained in its
divisor. Intersecting with $\E_2$ gives
$$0=<(a\E_1+{\cal H}).\E_2>=a-2$$
as by assumption, only $Q-P$ lies on $\E_2$ so $({\cal H}.\E_2)=-2$

A similar computation for $(\C_Q,f_Q)$ shows that 
$$\partial ((\C_Q,f_Q))=2\E'_2+{\cal H'}$$
where  here $C_Q=\E'_1+\E'_2$ are  translates of the $\E_1$ and $\E_2$
which split $C_P$ and ${\cal H}'$ is the closure of $div f_Q$.

The  cocycle  condition $div(f_p)+div(f_Q)=0$,  gives  ${\cal H}+{\cal
H}'=0$, as ${\cal H}$  and  ${\cal H}'$  are   the closures  of  these
divisors.      Adding      the     two     and      observing     that
$(\E_1-\E_2)-(\E_1'-\E_2')$ is homologous  to zero  in the fibre   (in
fact it is torsion) we get our result.

\end{proof}

\begin {rem}

This calculation works only locally as, a priori, it is not clear that
one can  find  a global decomposable  element  which has  the required
boundary. However, we will later give an argument  that this is indeed
the case.

\end {rem}

\subsection{Computing Signs}
 
In this   section we give a recipe    for computing the  signs  of the
boundary.  As shown  in the previous section, the  sign depends on the
position of  the 2-torsion points $O$  and $P-Q$ when the curve splits
into a product of  two elliptic curves.  The best way to describe this
is   to use  genus 2 theta    functions  with characteristics.  A good
reference for all the facts used here is \cite{kf}.

\bigskip

\noindent{\em Theta Functions}

\bigskip

The {\bf degree 2 Theta function with characteristic}
$[\epsilon,\epsilon']$ is defined to be
$$\Theta [\epsilon,\epsilon'] (z,\tau)=\sum_{N
\in  \ZZ^2}  exp (2\pi  i (\frac{1}{2}(N+\frac{\epsilon}{2}) \tau
(N+\frac{\epsilon}{2})^{tr}) + (N+\frac{\epsilon}{2})(z+\frac{\epsilon'}
{2})^{tr})$$
where $\tau=\pmatrix  \tau_1 & \tau_2 \\  \tau_2 & \tau_3 \endpmatrix$
is a point on  ${\frak  H}_2$, $z=[(z_1,z_2)]$  is a point  on  ${\Bbb
C}^2$  being thought  of as  a row  matrix and $[\epsilon,\epsilon']  =
[(\epsilon_1,\epsilon_2),(\epsilon'_1,\epsilon'_2)]$  is a  pair    of
points in $(\ZZ/2\ZZ)^2$.

To  a  characteristic    $[\epsilon,\epsilon']$   we associate     the
$2$-torsion point $I\frac{\epsilon'}{2}^{tr}+\tau
\frac{\epsilon}{2}^{tr}$ on the abelian surface $\A_{\tau}$ and call 
this the {\bf associated $2$-torsion point}.

A theta function with characteristic is called {\em odd} or {\em even}
depending  on  whether the corresponding function   is  an odd or even
function.             This         is         equivalent            to
$\epsilon_1\epsilon'_1+\epsilon_2\epsilon'_2$ being odd or even.

Up to  a     nowhere-vanishing  factor,  the   theta  function    with
characteristic is  the same  as the function  that  one would  get  by
translating   the theta function  $\Theta[(0,0)(0,0)](z,\tau)$  by the
associated $2$-torsion point. This shows that there  are six odd theta
functions and ten even theta functions.

The   functions        on   the     moduli        space    given    by
$$\Theta[(\epsilon_1,\epsilon_2),(\epsilon'_1,\epsilon'_2)](\tau)=
\Theta[(\epsilon_1,\epsilon_2),(\epsilon'_1,\epsilon'_2)](0,\tau)$$
are called the {\bf thetanullwerte}.  Since odd theta functions, being
odd  functions,  vanish   identically   at  $0$  there   are  no   odd
thetanullwerte.

Let $\Gamma_2=Sp_4(\ZZ)$ and $\Gamma_2(2)$ be the principal congruence
subgroup of level 2.   From the transformation  formula for  the theta
functions, one can see that the  thetanullwerte descend to give modular
forms  of weight  one-half on  ${\frak H}_2/\Gamma_2(2)$.   A well  known
theorem, see for   example \cite{vdg}, says that each  thetanullwert
vanishes precisely on  one  component  of  the moduli  of products  of
elliptic curves.  

The group  $\Gamma_2/\Gamma_2(2) \simeq Sp_4(\ZZ/2\ZZ)$  acts on these
theta functions or theta characteristics.  It acts transitively on the
set of  odd theta characteristics and can  be used \cite{ig} to get an
isomorphism   of  $Sp_4(\ZZ/2\ZZ)$    with  ${\frak    S}_6$   .   Let
$\Gamma_2(o)$  be the  inverse   image of  ${\frak  A}_6$ under   this
isomorphism and let $\omega \in Sp_4(\ZZ)$ be such that $\omega \notin
\Gamma_2(o)$ and $\omega^2=I$.

At a  point $\tau$ where  the polarization splits  into  a sum  of two
elliptic  curves,  $\C=\E_1+\E_2$, the existence   of the curve $\E_1$
implies that the period $\tau$ satisfies  a certain singular relation,
$v=(\alpha,\beta,\gamma,\delta,\epsilon)$. The        curve     $\E_2$
corresponds to  the relation $-v$. We can  make a  choice of $\E_1$ or
$\E_2$ and look at the  $\Gamma_2(o)$ translates, namely 
$$g\E_1g^{tr} \text{for all $g$ in } \Gamma_2(o).$$
This allows us to make a uniform choice of elliptic curve $\E_1$, with
respect to which  we can compute  the sign.  The sign does  not change
within components.

To compute the sign, we first need to make a choice of an embedding of
the genus 2 curve.  In Farkas and Kra \cite{kf} chapter VII, they make
a   choice  and with respect  to   their  choice, the   genus  2 curve
$\C_{P_1}$ is given by the divisor odd theta function
$$ \Theta[(0,1),(0,1)] (z,\tau)$$
and the six ramification points are  the two torsion points associated
to the characteristics
$$P_1=[(0,0),(0,0)] \,\,\,  P_2=[(1,0),(0,0)] \,\,\,  P_3=[(1,0),(1,1)]$$
$$P_4=[(1,1),(1,1)] \,\,\, P_5=[(1,1),(1,0)] \,\,\, P_6=[(0,0),(1,0)]$$

In the component of the moduli of  products of elliptic curves given
by $\tau=\pmatrix \tau_1 & 0 \\ 0 & \tau_3 \endpmatrix$, the odd theta
function splits as   

$$ \Theta[(0,1),(0,1)]((z_1,z_2),\tau)=\Theta
 [(0,0)] (z_1,\tau_1) \, \Theta [(1,1)] (z_2,\tau_3)$$
Since    the    zeroes  of    the   elliptic    theta  functions   are
$[(1,1)]=\frac{1}{2}(1+\tau_1)$  and $[(0,0)]=0$   respectively,   the
elliptic curves meet at the 2-torsion point
$$R=[(1,0),(1,0)]$$
From this one can see that the points $P_1,P_2$ and  $P_6$ lie one the
elliptic curve $z_1 \times[(0,0)]$ and the  points $P_3,P_4$ and $P_5$
lie on the elliptic curve $[(1,1)] \times z_2$.

So, for example,  if we choose  $P=P_1$ and  $Q=P_3$ and  consider the
element  of  $CH^2(\WW,1)$  given  by   $(\C_P,f_P)+(\C_Q,f_Q)$,   the
boundary of this element would have a sign of $-1$ on the component of
the  moduli where the curve   $\C_P$ splits in  to a  sum of two elliptic
curves meeting at the $2$-torsion point $[(1,0),(1,0)]$.

In   general,  to a component  (  which  corresponds to an  even theta
characteristic ) we  can associate a  $(3,3)$ configuration of the six
ramification points, namely    the first triple  corresponding  to the
points  lying on our  choice of  $\E_1$ and  the second to  the points
lying on $\E_2$.  To any pair of these  ramification points  we have a
Collino element of $CH^2(\WW,1)$ and the sign of the boundary on this
component would depend on the which triple $0$ and  $P-Q$ belong to. If
they belong  to the  same triple, then   the sign is  $0$ and  if they
belong to different  triples, the  sign  is $+1$ or  $-1$ depending on
whether $0$ or $P-Q$ lies on $\E_1$.

The sign   at  a different component  will  depend  on the  action  of
$\Gamma_2(o)$ on   the $(3,3)$ configuration   of the six ramification
points.   To  understand  this action,  we  observe that  this can  be
translated in to the well known \cite{ig},\cite  {vdgh} situation of the
action of  ${\frak A}_6$ on the  $(3,3)$ configuration of  the six odd
theta characteristics. This is  because the  points associated to  the
six odd theta charateristics are  nothing but the ramifications points
of the genus 2 curve given by divisor of  $$\Theta[(0,0),(0,0)](z,\tau)$$
As we used $$ \Theta[(0,1),(0,1)] (z,\tau)$$  to compute the action on
the six ramification points we have, we simply  translate our points by
the   point  $[(0,1),(0,1)]$ which  gives   us    the  six odd   theta
characteristics, and then use the description of the action on them.

The table  describes  the correspondence between  pairs  of  odd theta
characteristics    and   even  theta   characteristics.    We use  the
convention,  as in \cite{ig} that  the six odd theta chacteristics are
$$P_1=[(0,1),(0,1)]           \,\,\,     P_2=[(0,1),(1,1)]      \,\,\,
P_3=[(1,0),(1,0)]$$    $$P_4=[(1,0),(1,1)]  \,\,\,   P_5=[(1,1),(0,0)]
\,\,\, P_6=[(1,1),(1,0)]$$ With that convention, the class $(346,125)$
corresponds to  $[(1,1),(1,1)]$  The $g$    in  the leftmost    column
represents  the element of ${\frak A}_{6}$  which has to be applied to
$(346,125)$ to obtain the pair on the right. In \cite{ig} he makes the
isomophism  of ${\frak  S}_6$  with $Sp_4(\ZZ/2\ZZ)$  explicit,  which
along with some help from MAPLE, allows one to  do the computation and
make the following table.

\begin{table}
$$
\begin{tabular}{l|c|c|c}                                        
\hline

& g & Pair & Even Theta Characteristic  \\ \hline
I    & e        & (346,125) & 1111 \\
II   & (14)(26) & (123,456) & 1000 \\ 
III  & (56)(12) & (345,126) & 1100 \\  
IV   & (26)(34) & (124,356) & 0110 \\  
V    & (13)(56) & (145,236) & 0011 \\ 
VI   & (13)(45) & (156,234) & 1001 \\ 
VII  & (24)(56) & (235,146) & 0000 \\
VIII & (23)(34) & (246,135) & 0010 \\ 
IX   & (23)(45) & (256,134) & 0100 \\ 
X    & (14)(46) & (136,245) & 0001 \\

\hline
\end{tabular}
$$
\caption[]{Correspondence between (3,3) configurations and even Theta Characteristics}
\end{table}

Define the function

$$
\epsilon_{P,Q}(\tau ) = \begin{cases}
 -1 &\text  {if  $(O)$ lies on $(\E_1)_{\tau}$    and $(P-Q)$ lies  on
 $(\E_2)_{\tau}$} \\ 0 &\text  {if $(O)$  lies on $(\E_i)_{\tau}$  and
 $(P-Q)$ lies on $(\E_i)_{\tau}$ for $i=1,2$}\\ 1&\text {if $(O)$ lies
 on $(\E_2)_{\tau}$ and $(P-Q)$ lies on $(\E_1)_{\tau}$}
\end{cases} 
$$

From the  above discussion,  given  a  point $\tau$ on  the  component
corresponding         to     $[(1,1),(1,1)]$   we      can   determine
$\epsilon_{P,Q}(g\tau)$ for $g$ in $\Gamma_2(o)$.

\subsection{Computing the relation}

In this section, we restrict   the computations done in the   previous
sections to Shimura curves.

From  Hashimoto's \cite  {ha},  Thm 3.5,  explicit description  of the
Eichler   order  and Shimura curve   one  can explicitly  describe the
embedding  $\Psi:{\frak  H}   \hookrightarrow {\frak   H}_2$  which is
compatible with the  embedding of  the group  of units $\Gamma$   into
$\Gamma_2(o)$ :

\begin{thm}

The curve ${\bf X(D_0,N)}$ satisfies singular relations of the form 

$$ \tau_1+\tau_2+\-\frac{p-1}{4}\tau_3=0 $$
$$ 2aD\tau_2+(\tau_2^2-\tau_1\tau_3)+(a^2D-b)D=0 $$
giving rise to an embedding 

$$\Psi(\tau)=\frac{1}{p              \tau}\pmatrix              -{\bar
{\kappa}}^{2}+\frac{(p-1)aD}{2}\tau+D{\kappa}^2{\tau}^2              &
{\bar{\kappa}}-(p-1)aD\tau-D         \kappa       {\tau}^2          \\
{\bar{\kappa}}-(p-1)aD\tau-D\kappa{\tau}^2 & -1-2aD\tau+D{\tau}^2
\endpmatrix $$
where $\kappa=\frac{1+\sqrt p}{2}$

\end{thm}

From this one can see \cite {ha} that  the $\tau$ in the image satisfy
singular relations with invariants of the form $pX^2+4aDXY+4bDY^2$ for
$(X,Y)$ coprime.

Let $v$ be the singular relation given by the existence of an elliptic
curve   of degree $1$,  and let  $\Psi(\tau)$  be a point on $\Psi(\HH)$
which satisfies this relation. Then $\tau$ satisfies an equation given
by
$$\pmatrix  \Psi(\tau) & I_2 \endpmatrix  v \pmatrix \Psi(\tau) \\ I_2
\endpmatrix =0 $$
Since  as  functions   of   $\tau$, $\tau_1$,$\tau_2$,   $\tau_3$  and
$\tau_2^2-\tau_1\tau_3$  are  at  worst quadratic   when multiplied by
$p\tau$ and  since the  singular  relation between  the $\tau_i's$  is
linear, the   equation satisfied by $\tau$  is   quadratic.  Hence the
point $\tau$ is imaginary quadratic so is a CM point on $\HH$.

We  can look    at   the action  of $\Gamma_2(o)$ on $v$, namely
$$g \cdot v= g^{tr} v g \text { for $g$ in $\Gamma_2(o)$ } $$
The point $g\Psi(\tau)$, which is given by the action of $\Gamma_2(o)$
on $\HH_2$ need not  lie on $\Psi(\HH)$.   However, for some $g$ there
exists a point $\Psi(\tau')$ which satisfies
$$\pmatrix  \Psi(\tau')   &  I_2   \endpmatrix   g^{tr} v g    \pmatrix
 \Psi(\tau') \\ I_2 \endpmatrix =0 $$ 
If such a  point exists it will  be unique, and we  will  denote it by
$g(\tau)$. This has   the  same level   $2$ structure as   $g\tau$, so
$\epsilon_{P,Q}(g(\tau))=\epsilon_{P,Q}(g\tau)$.     Let   $\widetilde
{g(\tau)}$ be the point $\omega(g(\tau))$.  Then $\epsilon_{P,Q}(
\widetilde {g(\tau)})=-\epsilon_{P,Q}(g(\tau))$

The relation  we  get is  then  a signed sum  over  the $g(\tau)$, and
putting it all together gives us the following proposition.

\begin{prop}

Let   $(C_P,f_P)+(C_Q,f_Q)$  be   Collino's  element   in   the  group
$CH^2(\A_{\eta},1)$, where $\eta$ is the  generic  point of a  Shimura
curve $\X$ as before.  Let  $\tau$ be a  point on the moduli where the
curve $C_P$ splits as a sum of elliptic curves $\E_1+\E_2$, and make a
choice  of $\E_1$ as  before.  Assume further  that $O$ lies on $\E_1$
and $Q-P$ lies on $\E_2$. Then one has a relation in $CH^2(\W) \otimes
\Q$ of the form
$$\sum_{g      \in            \Gamma_2(o)}     \epsilon_{P,Q}(g(\tau))
(2(\E_2-\E_1)|_{g(\tau)}-2(\E_2-\E_1)|_{\widetilde {g(\tau)}}) \equiv 0$$
up to  the   relations  coming from  decomposable   elements and cycles
homologous  to      zero     in   the     special    fibres,     where
$\epsilon_{P,Q}(g(\tau))$ can be computed from the table in section $4.3$.

\end{prop}

\subsection{Isogenies}

In this section we  will explain how we  can use generic  isogenies to
get more relations between some cycles in $CH^2(\W)$.

Let  $g$ be an element  of $Sp_4(\Q)$. Let $\tau$   be a point. 
Then one has an induced isogeny
$$\Phi_g:\A_{g\tau} \rightarrow \A_{\tau}$$
of some degree $n^2$. For a fixed degree, there are only finitely many
classes of $g$ under the relation
$$g_1 \sim  g_2  \Longleftrightarrow \exists  \,  \gamma  \in Stab(\HH
\times \HH)|\,\,\, g_1=\gamma g_2$$

We can use this information to compute the relations  one would get by
applying such isogenies to  the original relation in Proposition $4.4$
. Equivalently,  we could apply  the isogenies to Collino's element to
get new elements and compute their boundaries.

In  the previous  section   we  have  used  the  embedding   $\Psi:\HH
\hookrightarrow \HH_2$. Let $g$ be in $Sp_4(\Q)$ of degree $n^2$ prime
to $2k$.    $g\Psi$  gives another  embedding  of  embedding   of $\HH
\hookrightarrow \HH_2$.    Let $\Gamma_g$ denote   the group  $g\Gamma
g^{-1} \cap Sp_4(\ZZ)$, where the action of $g$  is on the image of of
$\Gamma$ in $Sp_4(\ZZ)$. This acts on the new embedding of $\HH$.

Let $\X_{g}=\HH/\Gamma_g$ (or its compactification, if  necessary).
There is a map $\X_g \rightarrow \X$ given by $\tau \rightarrow g^{-1}
\tau$ and as above, a corresponding isogeny  $\Phi_g$ of the universal
abelian   surfaces over these   curves.   Let $\zeta$ denote Collino's
element in $CH^2(\A_g,1)$. We  can  compute relations resulting   from
this element in the universal family over $\X_g$, $\W_g$ and then push
the  relation down  using the   isogeny.  This will give   us some new
relations in $CH^2(\W)$.  A remark is  necessary when the curve $\X_g$
is not compact as  then it is not immediately  clear that  the isogeny
extends  to  a  morphism  over the   cusps.  However,  as   the cycles
supported in the  cuspidal fibres are torsion  as all these curves are
quotients of the upper half plane by congruence subgroups \cite{scl2},
one can  always  multiply the relation   by a suitable number  to kill
those cycles, and  then apply the  isogeny.  As  the isogeny preserves
the non-cuspidal points, there is no problem.

The relation that one obtains in  $CH^2(\W_g)$ are supported in fibres
where there is   an elliptic curve  of  degree $1$.   Let $v$   be the
singular  relation  representing   the elliptic  curve  of  degree  1,
$\E_1$. By translating by $\Gamma_2(o)$, we  can make a uniform choice
of the  elliptic curve  $\E_1$ over the  entire  moduli of products of
elliptic curves.  Let  $\tau$ be  a  point on  $\X_g$ satisfying  that
singular  relation.  Then, as described  above, the point $g^{-1}\tau$
on $\X$ satisfies   a  relation  corresponding  to  the  existance  of
elliptic curve of degree $n$ in the fibre. 

It is  possible for $g^{-1}\tau$  to lie  on the  image of some  other
point $\tau'$ as
$$g^{-1}\tau=\gamma g^{-1}\tau'  \text{\space   for some $\gamma$   in
$\Gamma(2k)$ } \Rightarrow g \gamma g^{-1}\tau=\tau'$$
However, we are looking at equivalence by the subgroup of finite index
in $g\Gamma g^{-1}$ given by $\Gamma_g$,  and these points need not be
$\Gamma_g$  equivalent. But since  $g$  gives rise   to an  isogeny of
degree coprime with $2k$ the signs of $\E_1-\E_2$ in all the points in
the fibre over a point on $\X$ are the same.

Conversely, a theorem of Kani  \cite{ka} asserts that if $\tau_0$ is
a point on $\X$ where there is an elliptic  curve of degree $n$ in its
fibre, then $\tau_0  = g^{-1}  \tau$ for  some point $\tau$  on $\X_g$
where there is  an elliptic curve of  degree  $1$ and $g$ a  primitive
isogeny of degree $n^2$.

Let $g_1,....,g_M$ be representatives  for the finitely many primitive
isogeny classes of  degree $n^2$. Then, since we  have made  a uniform
choice of $\E_1$ over the  entire moduli, we  have a uniform choice on
all the points $\tau$ on all the $\X_{g_i}$ which lie on the moduli of
products. Since the isogenies are of degree prime to $2k$, they do not
affect     the    level    structure,    and  so    one     can define
$\epsilon_{P,Q}(g^{-1}\tau)=\epsilon_{P,Q}(\tau)$  for  $g$  of degree
$n^2$.

Applying the corresponding   isogenies to the  relations on $\X_{g_i},
i=1,...M$ one has the following proposition, which is a generalization
of the proposion in the previous section.

\begin{prop}

Let $\tau_1$ be a  point  on $\X_{g_1}$ , where   $g_1$ is one of  the
representatives for the isogenies of  degree  $n^2$, and using  $\tau$
make a uniform  choice of elliptic curve $\E_1$  for the whole  moduli
space. Choose a point $\tau_i$ on each of the $\X_{g_i}$

Then, for each $i$,  there are relations in  $CH^2(\W) \otimes \Q $ of
the form
$$\sum_{g           \in      \Gamma_2(o)}               \epsilon_{P,Q}
(g_i^{-1}g(\tau_i))(2(\E_2-\E_1)|_{g_i^{-1}g(\tau_i)}                -
2(\E_2-\E_1)|_{\widetilde  {g_i^{-1}g(\tau_i)}}) \equiv 0$$
where $\E_1$ and $\E_2$  are now elliptic  curves of degree $n$ in the
fibre over $g_i^{-1}g(\tau)$.

\end{prop}    
 
\begin{rem}

While it may happen that $g_i \tau = g_j \tau'$ for some points $\tau$
and $\tau'$, this can  never happen on  an embedding of a  curve where
there is no generic  elliptic curve of degree  $n^2$ as  otherwise one
would have too many elliptic curves in the  fibre. This is because the
classes  of  elliptic  curves   are    linearly independent  in    the
N\'eron-Severi.  Both $g_i$ and  $g_j$ would  contribute  2 each.   If
generically  there  were  no elliptic   curves this  would  lead  to a
contradiction   as there   would   be  too many  linearly  independent
elements.

\end{rem}

Hence we get infinitely  many relations between some cycles  supported
in points where  there is an elliptic  curve of degree prime to  $2k$.
Since generically,  the  N\'eron-Severi is  of  rank three, when  this
happens the N\'eron-Severi jumps to rank 4, and the abelian surface is
necessarily isogenous to a  product  of isogenous CM elliptic  curves,
and there are CM cycles in those fibres.

In the next section we  will use the relations  above to get relations
between the CM cycles.

\section{Relations between CM cycles}

\subsection{Rewriting the relation}

In  this section we first describe  the cycles $\E_1-\E_2$ in terms of
the   CM cycles  and  show that  we  can modify  Collino's  element by
suitably chosen  decomposable elements to  get a relation only between
them. Then we get a more explicit description of the points $g(\tau)$.
  
In all our computations we will work  modulo cycles homologous to 0 in
the fibre. From weight considerations,  one  can see that such  cycles
map to 0 in the second intermediate Jacobian, though  they need not be
0 in the Chow group. Since the CM cycles themselves are defined modulo
such cycles, this is not really a restriction. For simplicity, we will
only  work with the    case of Collino's   original element. 

If  we  fix the  functions   $f_P$ and $f_Q$,    then the boundary  of
Collino's element looks like
$$\sum_{\tau} {\cal D}_{\tau}$$
where ${\cal D}_{\tau}$ is some element of the N\'eron-Severi group of
the fibre  over $\tau$.   $D_{\tau}$  can be written  in  terms of the
basis  for the rational  N\'eron-Severi given by $e_1$,$e_2$,$e_3$ and
$\Z_{\tau}$, 
$${\cal  D}_{\tau} = b_{\tau}^1e_1  +  b_{\tau}^2e_2 + b_{\tau}^3e_3 +
c_{\tau}\Z_{\tau}$$
where $\Z_{\tau}$  denotes   the CM cycle  if   the point
${\tau}$ is a CM point and is $0$ otherwise

\begin{lem}

  Let $\W$ be (the compactification of) the  universal family of abelian
  surfaces as before and let ${\cal W}$  be the generic fibre.  Recall
  that $e_1$,$e_2$ and $e_3$ generate the generic N\'eron-Severi. Then
  there are decomposable  elements $(e_i,f_{e_i})$  for $i=1,2,3$ such
  that
$$ \partial    (\sum_{i}      (e_i,f_{e_i}))=\sum_{{\tau}}   (\sum_{i}
b_{\tau}^{i}e_i).$$
Hence there is a relation of the form 
$$\sum_{\tau}  c_{\tau}\Z_{\tau}+\{hom\} \equiv 0$$
in $CH^2(\W)\otimes{\Bbb Q}$, where   $\{hom\}$ denotes the   image of
cycles homologous to zero in the fibre.

\end{lem}

\begin{proof} 
Let $\pi_{i}:\overline{e_i} \rightarrow  \X$  be  the maps from    the
closure of the   $e_i$ to the base $\X$.    Let $M_i$ be   divisors in
$CH^1(\W)\otimes {\Bbb Q}$ be  such that $(M_i.e_j)=\delta_{ij}$. Such
$M_i$ exist as   one can always  find  an  orthonormal basis  for  the
rational N\'eron-Severi.  Intersecting $M_i$ with $\sum_{{\tau}} {\cal
D}_{\tau}$ gives a relation in $CH^2(\overline{e_i})$ of the form
$$\sum_{\tau} b_{\tau}^i (\overline {e_i} \cap M_i)|_{{\tau}}$$
The direct   image   $$(\pi_i)_*(\sum_{\tau}  (\overline {e_i}    \cap
M_i)|_{{\tau}})=\sum b_{\tau}^i {\tau}$$
is a  rational equivalence of  points on $X$ so  is the divisor  of a
function, $f_{e_i}$.

These functions $f_{e_i}$ combined  with  the elements $e_i$ give  the
required   decomposable elements.  Subtracting  these   elements  from
Collino's element gives a relation of the form

$$\sum_{\tau} c_{\tau} \Z_{\tau} +\{hom\} \equiv 0$$

\end{proof}

\begin{rem}
  
At the  cuspidal points, the cycles that  remain after subtracting off
the  boundary  of these decomposable  elements   are orthogonal to the
closure of the  generic N\'eron-Severi  and  as mentioned before   are
torsion in the Chow group of the special fibre  itself, so do not make
a contribution after tensoring with the rationals.

One  could also say things more  intrinsically by  using the fact that
there  is a projector  in the  ring  of correspondences  with rational
coefficients on ${\cal A}$ which takes the varation of Hodge structure
$R^2\pi_{*}({\Bbb Q})$ onto  $Sym^2(R^1p_{*}({\Bbb Q}))$  and the fact
that the CM cycles lie in this part in  the fibres over the CM points.
We could then apply this correspondence to the element of $CH^2(\A,1)$
itself to get a relation involving only CM cycles.

\end{rem}

As a result  of this lemma  one  sees that  since one gets  a relation
between CM cycles which all  lie  in $CH^2_{hom}({\bf  W})$, we get  a
relation there and not merely in the Chow group.

We would     now  like to   determine the    points  ${\tau}$  and the
coefficients $c_{\tau}$. To determine  $c_{\tau}$ we have to write the
CM cycle in terms of the  basis for the primitive N\'eron-Severi given
by          ${\bar        e_i}          ,i=2,3$    and          ${\bar
e_4}=\frac{1}{2}(2e_4-<e_4.e_1>e_1)$, where    $e_4=\E_1$.       Since
$<e_4,e_4>=0$, one has  $({\bar  e_4}.{\bar e_4})=1$.  Define $r$  and
$s$  by    $$({\bar e_4}.{\bar e_2})=r   \,\,\,\,\,  ({\bar e_4}.{\bar
e_3})=s$$
so the intersection matrix looks like

$$\begin{array}{c|ccc}

      & {\bar e_2}   & {\bar e_3}   & {\bar e_4}\\ \hline
     {\bar e_2} & p  & 2aD & r\\
     {\bar e_3} & 2aD & 4bD &s \\
     {\bar e_4} & r            & s   & 1

\end{array}
$$

The idea now is to compute the  CM cycle up  to a multiple using
the property that  it is  orthogonal  to the generic  N\'eron-Severi. 
Define $d_{\tau}$  to   be the  smallest multiple   of  the  CM  cycle
$\Z_{\tau}$   which  lies  in     the  $\ZZ$-span   of  ${\bar   e_i},
i=2,3,4$. Then one has

\begin{prop}

Let $d_{\tau}$ be as above, then 

$$d_{\tau}\Z_{\tau}       =       \frac{aDs-2bDr}        {2}     {\bar
e_2}+\frac{2aDr-ps}{2}{\bar e_3}+2D{\bar e_4}$$
up to a sign.

\end{prop}

\begin{proof}
Let 
$$\Z_{{\tau}}=x_2{\bar e_2}+x_3{\bar e_3}+x_4{\bar e_4}$$
Then $(\Z_{\tau},{\bar e_i})=0$ for $i=2,3$, so we get 2 equations
$$0=px_2+2aDx_3+rx_4$$
$$0=2aDx_2+4bDx_3+sx_4$$
Eliminating the variable $x_2$ gives 
$$x_3=\frac{2aDr-ps}{4D}x_4$$
and similarly

$$x_2=\frac{2aDs-4bDr}{4D}x_4$$
so the vector 
$$(\frac{2aDs-4bDr}{4D},\frac{2aDr-ps}{4D},1)$$
is orthogonal to the  generic N\'eron-Severi hence a rational multiple
of the CM cycle.

The smallest multiple of the vector such that it is an integral linear
combination will  be the $d_{\tau}$ times the  CM cycle or half the CM
cycle depending  on the discriminant, at  least up to  a sign.  Define
$b_{r,s}$ as the  multiple which  gives it,   hence it is  either  the
smallest multiple or twice the smallest multiple.

From the description  above, and the fact  that $s$ is even,  which we
will   see    presently,   it   is   clear    that   $b_{r,s}|4D$  and
$c_{\tau}=c_{r,s}=4D/b_{r,s}$.  The smallest multiple is determined by
$s_1=p(s/2)-aDr$, if $s_1$ is even,  then it  is  $D$ otherwise it  is
$2D$.

The proposition will then follow from the following lemma

\begin{lem}

$|b_{r,s}|=2D$ always.

\end{lem}

\begin{proof}

Let $v=(\alpha,\beta,\gamma,\delta,\epsilon)$ be the singular relation
corresponding  to the   elliptic   curve  $\E_1$. From   Hashimoto's
explicit description of the embedding, one can compute  $r$ and $s$ to
be
$$r=\beta-2\gamma-\frac{1-p}{2}\alpha$$
$$s=2aD\beta-2\epsilon-2D\delta(a^2D-b)$$
Let 
$$ s_1=p(s/2)-aDr =( \frac{(p-1)}{2}aD\alpha-(p-1)aD\beta - 2aD\gamma
+ ((p-1)a^2D)\delta + p\epsilon)$$
It    turns out that $\tau$   satisfies   an equation of  discriminant
$(1/p)(s_1^2-4D(\frac{p-r^2}{4}))$. This is even if  and only if $s_1$
is even.  Hence  $b_{r,s}=2D$ always, as  if $s_1$  is even, then  the
smallest multiple is  $D$  hence $b_{r,s}=2D$  is twice the   smallest
multiple.   If $s_1$ is  odd,then  the smallest  multiple is  $2D$ and
$b_{r,s}$  is  the  smallest   multiple.  This  also  says  that   the
discriminant of  $\tau$ is $(1/p)(s_1^2-4D(\frac{p-r^2}{4}))$ up  to a
square factor.

\end{proof}   

This concludes the proof of the proposition. 

\end{proof}

Armed with this information we can compute discriminants of the orders
that appear in  the relation. They are  the points on the moduli where
the N\'eron-Severi  is generated by the  generic  elements along with an
elliptic curve of degree $1$.

\begin{prop}

The discriminants of the points $\tau$ where there  is a elliptic curve
of degree $1$ are 
$$\Delta/d^2=\frac{P(s,r)-4D}{4d^2}=(p(s/2)^2-aDrs+bDr^2-D)/d^2$$
where $r$ and $s$, where $r$ and $s$ run  through all possible numbers
satisfying $\frac{p-r^2}{4}  \in \ZZ_{\geq 0}$, $\frac{4bD-s^2}{4} \in
\ZZ_{\geq 0}$ and $d$ is such that $d^2|\frac{P(s,r)-4D}{4}$ 
\end{prop}

\begin{proof}

The  idea is to compare  the determinants of the intersection matrices
of the  two bases  of  the rational  N\'eron-Severi.  Using  the  basis
coming    from ${\bar    e_1},{\bar      e_2}$  and  the   CM    cycle
$d_{\tau}\Z_{\tau}$, at a  point  of discriminant $\Delta_0$  one gets
the determinant of the intersection matrix to be

$$-16D^2d_{\tau}^2\Delta_0$$  
On the other  hand, from the computation  using the  other basis given
by ${\bar e_1},{\bar  e_2},{\bar e_3}$ and  the fact that  the change
of   basis   matrix  has   determinant  $b_{r,s}^2=4D^2$,  a   simple
calculation shows that the determinant is  $$(4D-P(s,r))4D^2$$
where $P(s,r)$ is the quadratic form $ps^2-4aDrs+4bDr^2$.

Comparing the two shows that 
$$\Delta_0=\frac{P(s,r)-4D}{4d_{\tau}^2}=\Delta/d_{\tau}^2$$
where $\Delta=\frac{P(s,r)-4D}{4}$

If  $v_1$ and $v_2$  are two  elements  of the primitive N\'eron-Severi,
then the lattice  generated by them is  the N\'eron-Severi lattice of an
abelian    surface  with  multiplication   by   an   Eichler order  of
discriminant  $\frac{v_1^2v_2^2-(v_1.v_2)^2}{4}$.   This  explains why
$\frac{p-r^2}{4}$ and $\frac{4bD-s^2}{4}$ are in $\ZZ$.

Conversely,  if $r$ and  $s$  satisfy the conditions that $\frac{p-r^2}{4}$
and   $\frac{4bD-s^2}{4}$     are in    $\ZZ$  and    $d^2$    divides
$\frac{P(s,r)-4Dn^2}{4}$ then one can see that the element
$$\E=\frac{1}{2D}(d\Z_{\tau}-\frac{aDs-2bDr}{2}{\bar
e_2}-\frac{2aDr-ps}{2}{\bar e_3})$$
where   $\tau$     satisfies      an  equation    of     discriminant
$\frac{P(s,r)-4D}{4d^2}$,  is   a   primitive  integral  element  of
invariant  $n^2$,  hence, from  \cite{ka}, one   sees  that it   is an
elliptic curve of degree $1$.

\end{proof}

Since $P(s,r)$ is positive definite, the  discriminant is negative for
only finitely many values of $s$ and $r$.

To generalize Proposition $5.5$ for the relations we get from applying
isogenies, we have to carry out the same computation using an elliptic
curve  of odd degree  $n$   instead of $1$,  and   one gets that   the
discriminants that appear are
$$\Delta=\frac{P(s,r)-4Dn^2}{4d^2}=\frac{p(s/2)^2-aDrs+bDr^2-Dn^2}{d^2}$$

\subsection{Heegner cycles}

In this section we define the {\bf Heegner cycles}  which are sums of
certain CM cycles with  the same discriminant.   As it turns  out, our
final  result can be expressed in  terms of these  cycles. A reference
for the facts used here is \cite{bd} or \cite{gz}.

Let  $\OO$  be,  as  before, the  Eichler order   in the  quaternionic
division   algebra  (or  $M_2({\Bbb Q}$).   Let  $K$   be an imaginary
quadratic field  and let  $\OO_{\Delta}$ be  an order  of discriminant
$\Delta$. If $\tau$ is a CM point, then $\tau$ determines an embedding
of $K \hookrightarrow \B$ and $\tau$ is called a {\bf Heegner point of
discriminant $\Delta$} if $$K \cap \OO=\OO_{\Delta}$$

Composing with the orientation on the Eichler order ${\frak o}_{\ell}$
gives an  {\em orientation} on the  Heegner point, namely a surjective
linear map
$$\kappa_{\ell}:\OO_{\Delta} \rightarrow {\Bbb F}_{\ell}$$
if $\ell | N$ or 
$$\kappa_{\ell}:\OO_{\Delta} \rightarrow {\Bbb F}_{\ell^2}$$
if $\ell | D_0$.

For a given orientation and level structure, it is well known
\cite{gz} that there are $h_{\Delta}$ Heegner points  of discriminant 
$\Delta$, where  $h_{\Delta}$  is the   class  number of  the   order.
Further, the  Class group $Pic(\OO_{\Delta})$     acts on the   points
preserving the orientation. The  Atkin-Lehner operator $w_{\ell}$, for
$\ell | D$ also acts on the Heegner points by flipping the orientation
at $\ell$ and also changing the ideal class.   The action of the group
$\Omega_{D} \ltimes  Pic(\OO_{\Delta})$ is transitive showing that for
a fixed level there are $2^{t}h_{\Delta}$ Heegner points.
 
Let $m$ be a number  such that 
$$m \equiv \begin{cases} a\,\,square\,\,\, mod \,\, \ell \,\,\,\,  if
\,\,\, \ell |N  \\   a\,\, nonsquare\,\,\,  mod \,\, \ell  \,\,\,\,  if
\,\,\, \ell |D_0 \end{cases}$$
and let   $\Q_{\ell,m}$ be the  extension  of $\Q_{\ell}$  obtained by
adding the  roots of the equation $x^2-m=0$.   For $\ell | N$  this is
simply $\Q_{\ell}$ and for $\ell | D_0$ this is a quadratic extension.
Let      $\ZZ_{l,m}$     be  the        ring     of    integers    and
$\OO_{\Delta,\ell,m}=\OO_{\Delta} \otimes \ZZ_{l,m}$.  Then
$$
\OO_{\Delta,\ell,m} / \ell \OO_{\Delta,\ell,m} = 
\begin{cases}{\Bbb  F}_{\ell} \oplus {\Bbb  F}_{\ell} & \text {if $\ell | N$}
\\ {\Bbb F}_{\ell^2} \oplus {\Bbb F}_{\ell^2} & \text {if $\ell | D_0$}
\end{cases}$$
and  the different orientations correspond  to the different canonical
maps or  equivalently   the different  primes  lying  over  the primes
dividing $D$.

From  that one can see that  the orientation classes  are in bijection
with solutions $\mu$ mod $2D$ of the equation
$$\mu^2 \equiv m\Delta  \,\,\,mod \,\,4D$$
as $\mu$ will give a solution of 
$$ \mu^2 \equiv m\Delta  \,\,\,mod \,\,4\ell$$
Since $m^{-1}$ is a square  in $\ZZ_{\ell,m}$, $m^{-1}=t_{\ell}^2$ for
some $t_{\ell} \in \ZZ_{\ell,m}$, so the ideal
$${\frak L}=(\ell,\frac{\mu t_{\ell}+\sqrt{\Delta}}{2})$$
is a  prime lying over $\ell$  in $\ZZ_{\ell,m}$. Conversely, given an
oriented Heegner point, the  kernels of  the various orientation  maps
give   rise  to the  different   primes and hence    a solution of the
equation.

The  choice of $m$ is  not that important.  Any  number satisfying the
conditions satisfied by $m$ will give the same result.

A Heegner point  on  $\X$ also comes with   the data of  a  level $2k$
structure, and this has a projection onto the level 2 structure $\LL$,
and we can further consider the  class of $\LL$,  $[\LL]$ which is the
even theta characteristic  corresponding to it.  We can make a uniform
choice of an elliptic curves $\E_1$ and $\E_2$  of degree $1$, or more
generally of degree  $n$, for $n$ odd,  over each of  these classes. A
{\bf Heegner  Cycle}   is the CM  cycle    in the fibre over  such   a
point. For a given level  $2$ structure and  orientation there are two
Heegner cycles corresponding  to which of the  two pairs of  three $2$
torsion points lies on   $\E_1$.  We define   the  {\em sign}  of  the
Heegner cycle  to  be  positive if  the   component of $\E_2$  in  the
direction of   the  CM  cycle   is  positive.   This   depends on  the
configuration of the $2$-torsion points.

From class  field theory, one   knows that  $Pic(\OO_{\Delta})  \simeq
Gal(H_{\Delta}/K)$, where $H_{\Delta}$ is  the   ring class field   of
$\OO_{\Delta}$. From this  one can  realise the  action of $Pic$  as a
Galois action, and one  can see that  the action on the Heegner cycles
does   not change the   sign.  The  action   of the  Fricke involution
$w_D=\prod_{\ell} w_{\ell}$ is,  up to an element  of $Pic$,  given by
complex conjugation and it flips the orientation and the configuration
of the pair of three $2$-torsion points and hence the sign. 

We   then can define a {\bf   Heegner cycle  of discriminant $\Delta$,
level $[{\cal L}]$ and orientation $\mu$} to be a sum 
$$\Z_{\Delta,[\LL],\mu,-\mu}    = \sum_{\LL \in [\LL]}     \sum_{{\frak a}
    \in Pic(\OO_{\Delta})}  \Z_{\frak a,\LL,\mu,I}  - \sum_{{\frak  a}   
\in Pic(\OO_{\Delta})} \Z_{\frak a,\LL,-\mu,II}$$
where the $I$ and $II$ denote the possible configurations.  This cycle
is invariant  under  the action of $Pic$   and changes sign  under the
action   of      the        Fricke    involution so
$$\Z_{\Delta,[\LL],\mu,-\mu}= \sum_{\sigma \in  Gal(H/\Q)}  \Z_{{\frak
a}_0,\LL,\mu_0,I_0}^{\sigma}$$

The  cycle $\Z_{\Delta,[\LL],-\mu,\mu}$  is   a different cycle  which
corresponds  to the same  configuration but  opposite orientation.  To
reduce  notation, and since both  the  configurations  occur, we  will
suppress the $I$ and $II$.

For  full level $2k$  structure we define  the Heegner cycle to be the
sum  over   all  the points  over   the  point determined by  $({\frak
a},\LL,\mu)$  of the Heegner cycles, and  we use the  same notation to
denote it.

In our situation,  if $r$ and $s$  denote the intersection  numbers of
$\E_1$ with  ${\bar e_2}$ and  ${\bar e_3}$,  then  they determine the
orientation of the Heegner  cycle  corresponding to  it, which is   of
discriminant $\Delta=\frac{ps^2-4adrs+4bDr^2-4D}{4d^2}$ for some $d$.

The prime $p$ is a number  like $m$ above, namely  a square mod $\ell$
for  $\ell | N$   and a  non-square mod  $\ell$   for $\ell |D_0$.  By
`completing the  square'  one can    then see that    $s_1=p(s/2)-aDr$
satisfies an equation of the form
$$x^2=p\Delta d^2 \,\,\,mod \,\,4D$$
and hence determines the  orientation !

\section{Main Result}

In this  section we state our main  result.  We  will assume that we
always have level $2k$ structure for some  odd integer $k$ coprime to
$N$.
  
\begin{thm}

Let ${\bf X(D_0,N)}$ be a Shimura curve parametrising abelian surfaces
with endomorphism ring  an Eichler  order  of level $N$ in  a division
algebra   of discriminant  $D_0$.  Let  ${\bf   W(D_0,N)}$ denote (the
non-singular compactification of)  the universal abelian  surface over
${\bf   X(D_0,N)}$.   Let $p,a,b$ be   the  invariants which determine
Hashimoto's  model  and   let  $P$  and  $Q$  denote  two  $2$-torsion
points. 

Then there are relations in $CH_{hom}^2({\bf W(D_0,N)})\otimes
\Q$   of the  form
$$  \sum_{[\LL]}  \epsilon_{P,Q}([\LL]) \sum_{r,s}
\sum_{d^2|\Delta} d \Z_{\Delta/d^2,[\LL],s_1/d,-s_1/d} \equiv 0 $$
where  $\LL$ runs   through   all  even  theta   characteristics   and
$\epsilon_{P,Q}([\LL])$ is the sign function as in Section $4.3$, $r,s
\in   \ZZ$  such    that  $\frac{n^2p-r^2}{4}   \in   \ZZ_{\geq   0}$,
$\frac{4bDn^2-s^2}{4} \in \ZZ_{\geq  0}$, $n$ is coprime   with 
$2k$, $r,s,n$ mutually  coprime, $s_1=p(s/2)-aDr)$ and
$$\Delta=\frac{ps^2-4aDrs+4bDr^2-4Dn^2}{4}$$
with  $d$ running through all $d^2|\Delta$   such that $\Delta/d^2$ is
still a discriminant.

\end{thm}

\begin{proof} 

This is just a consequence of  putting together the statements in the
previous sections.

\end{proof}

\begin{rem}

For different  choices of the  prime $p$ and numbers   $a$ and $b$ one
could get possibly different  relations.  One way of  interpreting the
Hashimoto model is  by  observing that it gives   an embedding of  the
Shimura curve in to the Siegel upper half space.  Different choices of
$p$,$a$ and $b$ correspond to  possibly different embeddings and could
result in different relations.

\end{rem}

\subsection{Examples}

I. Suppose  $D_0=2.3$,  $N=1$ so  $D=D_0N=6$  and $n=1$.   The  triple
$(p,a,b)=(5,2,5)$ determines   an embedding.  The  possible values for
$(r,s)$ are $$\{ (1,0) , (1,2), (1,4) , (1,6) , (1,8), (1,10) \}$$
Out of these only $(1,4)$ and $(1,6)$ give rise to negative values for
$\Delta$, which are $-4$ and $-3$ respectively.

Suppose the two torsion points are the ones associated to 
$$P = [(0,0),(0,0)] \text  { and } Q = [(1,0),(0,0)]$$
so they  correspond  to  the odd characteristics   $[(0,1),(0,1)]$ and
$[(1,1),(1,0)]$ respectively. Then from the table one can read off the
signs at the various components of the moduli of products corrsponding
to the  even theta  characteristics.   In this case,  for example, the
signs  are $-1$  at $[(1,1),(1,1)],[(0,1)(0,0)]$ and  $[(0,0),(1,0)]$,
$1$  at  $[(1,0),(0,0)],[(0,1),(1,0)]$  and $[(0,0),(1,1)]$  and   $0$
elsewhere.

So one has a relation in $CH^2( \bf {W(2,3)}\otimes \Q )$ of the form
$$\Z_{3,[1000],3,-3}+\Z_{4,[1000],2,-2}
+\Z_{3,[0110],3,-3}+\Z_{4,[0110],2,-2} $$
$$+\Z_{3,[0011],3,-3}+\Z_{4,[0011],2,-2} 
-\Z_{3,[1111],3,-3}- \Z_{4,[1111],2,-2}$$
$$-\Z_{3,[0100],3,-3}- \Z_{4,[0100],2,-2} 
-\Z_{3,[0010],3,-3} - \Z_{4,[0010],2,-2} $$
$$+\Z_{3,[1000],-3,3}+\Z_{4,[1000],-2,2} 
+ \Z_{3,[0110],-3,3}+\Z_{4,[0110],-2,2} $$
$$ +\Z_{3,[0011],-3,3} + \Z_{4,[0011],-2,2} 
-\Z_{3,[1111],-3,3}- \Z_{4,[1111],-2,2}  $$
$$- \Z_{3,[0100],-3,3} - \Z_{4,[0100],-2,2} 
- \Z_{3,[0010],-3,3} - \Z_{4,[0010],-2,2} =0$$
where for $[abcd]$ denotes the characteristic $[(a,b),(c,d)]$.

\noindent II. Different choices of  $(p,a,b)$ leads to different 
embeddings of   the  modular curves   and hence   new relations.   For
example, if  $D=26$, $(p,a,b)  =  (5,2,21)$, then there  are relations
between     the $\Delta    's$    given by  $(r,s,\Delta)=(1,18,-11)$,
$(1,20,-20)$ and $(1,22,-19)$.  However, if $(p,a,b)$ is $(149,19,63)$
then one has relations between $(r,s,\Delta)$ of the form $(1,6,-11)$,
$(3,20,-24)$, $(7,46,-11)$ and $(9,60,-8)$. In particular, one may get
a  relation    involving more than   one  Heegner  point of   the same
discriminant.

\begin{verse} 
Ramesh Sreekantan\\
Department of Mathematics\\
Duke University\\
Durham, NC 27708\\
{\bf email:}ramesh@math.duke.edu\\
\end{verse}

\end{document}